\documentclass{article}%
\usepackage{amsmath}
\usepackage{amssymb}
\usepackage{geometry}
\usepackage{graphicx}
\usepackage{amsfonts}%
\providecommand{\U}[1]{\protect\rule{.1in}{.1in}}
\newtheorem{theorem}{Theorem}[section]

\newtheorem{corollary}[theorem]{Corollary}

\newtheorem{definition}[theorem]{Definition}

\newtheorem{lemma}[theorem]{Lemma}

\newtheorem{proposition}[theorem]{Proposition}

\newenvironment{proof}[1][Proof]{\noindent\textbf{#1.} }{\ \rule{0.5em}{0.5em}}
\begin{document}

\title{Quandle-like Structures From Groups}
\author{Sriram Nagaraj\\Department of Electrical Engineering\\The University of Texas, at Dallas\\Richardson, TX 75083\\sriram@student.utdallas.edu}
\date{}
\maketitle

\begin{abstract}
We give a general procedure to construct a certain class of "quandle-like"
structures from an arbitrary group. These structures, which we refer to as
pseudoquandles, possess two of the three defining properties of quandles. We
classify all pseudoquandles obtained from an arbitrary finitely generated
abelian group. We also define the notion of the kernel of an element of a
pseudoquandle and prove some algebraic properties of pseudoquandles via its kernels.

\noindent%
\bigskip

\noindent%

\end{abstract}

\section{{\protect\LARGE Introduction}}

Quandles are algebraic structures that have been successfully employed in the
study of knots and links since their three defining axioms correspond to the
three Reidemeister moves. The presentations of quandles obtained from
knot/link diagrams are defined by considering the arcs as generators and the
crossings of the arcs as the relations. Quandles were first introduced by
Joyce in \cite{J-1}, although the idea of a rack has been available since the
1950's, particularly through the work of Conway and Wraith (see \cite{PN} for
more details). As a general idea, racks and quandles are the structures
obtained from a group $G$, where the group operation is replaced by
conjugation (or n-fold conjugation). This has its own advantages from the
group-theoretic point of view since the new structures have (at least in the
finite case) nice combinatorial properties. There are a set of (co)homological
ideas for theory of racks/quandles as developed in \cite{PN} and \cite{CJKS}.
We refer the reader to \cite{N-1},\cite{N-2},\cite{NN} for further discussion
of the properties of quandles, which we use later in this paper.

In this introductory section, we recall some basic definitions and examples of
racks and quandles, although the canonical reference for the subject is the
work by Joyce in \cite{J-1}, \cite{J-2}. We also introduce the idea of
pseudoquandles as algebraic structures satisfying two of the three defining
axioms of a quandle and show that we can obtain commutative pseudoquandles
from an arbitrary group $G$ (thus, although pseudoquandles do not correspond
to all the Reidemeister moves, they are easily obtained from any group). In
the following section, we classify all pseudoquandle obtained from an
arbitrary finitely generated Abelian group using the fundamental theorem of
Abelian groups. We also introduce the idea of the pseudoquandle matrix, an
extension of the ideas introduced in \cite{NV}, and \cite{NH} for the case of
quandles. In the final section, we define the kernel $ker(p)$ and cokernel
$coker(p)$ of any element $p$ of a pseudoquandle $P$ and prove some properties
of $P$ via the the kernels $ker(p)$ for each $p$ in $P$. In essence, the idea
of kernels of a pseudoquandle is an abstraction of various concepts introduced
in \cite{N-2}, \cite{NV}, and \cite{NH}. Basic results proven there may be
generalized using the notion of kernels.

\begin{definition}
A \textbf{quandle} is an algebraic structure $Q$ with a closed binary
operation $\ast:Q\times Q\longrightarrow Q$ that satisfies the following three
axioms$:$

$\mathbf{(i)}$ For all $q\in\ Q,$
\[
q\ast\ q=q
\]

$\mathbf{(ii)}$ $Q$ is a self distributive $($from the right$)$ with $\ast
$\ as an operation, i.e. for all $p,q,r\in\ Q$%
\[
(p\ast q)\ast r=(p\ast r)\ast(q\ast r)
\]

$\mathbf{(iii)}$ For each $p,q\in\ Q$ there is a unique $r\in\ Q$ such that
$p=r\ast q$
\end{definition}

If $Q$ satisfies axiom $\mathbf{(ii)}$ and $\mathbf{(iii)}$ above, we call it
a \textbf{rack}. We also define $Q$ to be a \textbf{pseudoquandle} if $Q$
satisfies $\mathbf{(i)}$ and $\mathbf{(ii)}$\textit{.} Note that axiom
$\mathbf{(iii)}$ above is equivalent to the following: For each $q\in\ Q$, the
map $\ast_{q}$\ from $Q$ into itself defined by:
\[
\ast_{q}(p)=p\ast q
\]
is bijective. Hence, in view of this, we may define an inverse (or dual)
operation $\ast^{-1}$ of $\ast$ which satisfies
\[
(p\ast q)\ast^{-1}q=p
\]
for each $p,q\in Q.$ We now give some standard examples of quandles. Most of
these examples can be found in quandle related literature.

The most basic quandle is the trivial quandle obtained from any set $Q$ with
operation
\[
a\ast b=a
\]
for all $a,b\in\ Q$. It can easily be seen that $Q$ a quandle. If $|Q|=n$, it
is called the trivial quandle of order $n$ denoted by $T_{n}$.

The (classical) example of a quandle is the quandle obtained from\textbf{\ }an
arbitrary group $G$. By fixing an integer $n$ and defining the operation
$\ast$ on $G$ as:%
\[
g\ast h=h^{-n}gh^{n}%
\]
we can see that $G$\ is a quandle with $\ast$\ as the quandle product. $($This
is $n$-fold conjugation in the group $G)$

There is also a class of quandles called Alexander quandles. An Alexander
quandle is a module $A$ over the ring $%
\mathbb{Z}
\lbrack t,t^{-1}]$ of formal Laurent polynomials with quandle product $\ast$
given by:%

\[
a\ast b=ta+(1-t)b\text{ for all }a,b\in\ A
\]

The inverse operation $\ast^{-1}$ in this case is given by:%

\[
a\ast^{-1}b=t^{-1}a+(1-t^{-1})b\text{ for all }a,b\in\ A
\]

For more details on the classification of finite Alexander quandles, the
reader can refer to \cite{N-1}. In the same vein as Alexander quandles, we can
define symplectic quandles as follows:

Let $M$ be a module over a ring $R$ $($with characteristic $\neq$ 2$)$. Let
$\langle,\rangle$ be an anti-symmetric bilinear form $\langle,\rangle:M\times
M\rightarrow R$ from $M\times M$ into $R$. Define the $\ast$ operation as
follows:%
\[
x\ast y=x+\langle x,y\rangle\text{ }y
\]
for all $x,y\in\ M$. It can be shown that $M$ with $\ast$ as a binary product
is a quandle called a symplectic quandle (see \cite{NN} for more
details\medskip).

As a final example, Let $Q=%
\mathbb{Z}
/n%
\mathbb{Z}
$. By defining $\ast$ on $Q$ as $i\ast j=2j-i$ (mod $n$) for all $i,j$ $\in$
$Q$, $Q$ becomes a quandle called a dihedral quandle of order $n$. Homological
methods as applied to these quandles have been developed in \cite{PN}

\bigskip

Now, we concentrate on the structures of our study, pseudoquandles obtained
from an arbitrary group $G$. Before we begin, we fix some notation used for
the rest of the paper.

\bigskip

\textbf{Notation}: Let $G$ be any group and $H$ a normal subgroup of $G$.
Define $G^{\vartriangleleft}$ as follows:%

\[
G^{\vartriangleleft}=\{H\text{ }|\text{ }H\text{ }\unlhd\text{ }G\}
\]
We label elements of $G^{\vartriangleleft}$ by alphabets $x,y,z...$etc. even
though they are actually normal subgroups. Let us define an operation $\ast$
on $G^{\vartriangleleft}$ as follows:%
\[
x\ast y=\{ab|\text{ }a\ \in x,b\ \in\ y\}
\]

This operation is simply the multiplication of two (normal) subgroups of $G$.
We then have the following:

\begin{proposition}
The set $G^{\vartriangleleft}$ is a commutative monoid with operation $\ast$
as defined above. Moreover, this operation is self distributive and every
element in $G^{\vartriangleleft}$ is idempotent with with respect to the
$\ast$ operation. Hence, $G^{\vartriangleleft}$ is a commutative pseudoquandle
with $\ast$ as the product.
\end{proposition}

%

\begin{proof}%
Let $x,y$ $\in\ G^{\vartriangleleft}$. Since all of the elements of
$G^{\vartriangleleft}$ are normal subgroups, their $\ast$ product is also a
normal subgroup, and this operation is commutative by definition of normal
subgroups. Since for every (normal) subgroup $g$ of $G$, $g\ast g=g$, every
element of $G^{\vartriangleleft}$ is idempotent. For self distributivity, we
have the following string of equalities:
\[
(y\ast x)\ast(z\ast x)=y\ast x\ast z\ast x\text{ }(\text{by assocativity of
the }\ast\text{ product obtained from }G)
\]%
\[
=y\ast z\ast x\ast x\text{ }(\text{since }x,y,z\text{ are normal})
\]

\[
=y\ast z\ast x=(y\ast z)\ast x
\]
for all $x,y,z$ $\in$ $G^{\vartriangleleft}.$%
\end{proof}%
\medskip

Whenever we refer to the \textbf{pseudoquandle obtained from a group }$G$, we
mean $G^{\vartriangleleft}$ with product as given above\footnote{Note that
$G^{\lhd}$ is self-distributive from both the right and left due to
commutativity}. We shall henceforth denote this structure as $P_{G}$. Although
$P_{G}$ constructed from any group $G$ is a pseudoquandle, even for the most
basic groups, $P_{G}$ is not a true quandle. For instance, if $G=\{\pm1,\pm
i,\pm j,\pm k\}$, the quaternion group, evey subgroup is normal and hence,
\[
P_{G}=\{\{1\},\{\pm1\},\{\pm1,\pm i\},\{\pm1,\pm j\},\{\pm1,\pm k\},G\}
\]
Consider $p=\{\pm1,\pm i\}$ and $q=\{\pm1,\pm j\}$; we see that there is no
(normal) subgroup $r$\ in $P_{G}$ such that $p=r\ast q$.

Since normal subgroups are those objects in the group which are invariant
under conjugation, we see that this construction resembles the classical
ideology of studying quandles, namely obtaining quandles by conjugation in the group.

\section{Pseudoquandles obtained via finitely generated abelian groups}

We begin this section with the following observations which motivates the
proof of the main result regarding the classification of $P_{G}$ for any
finitely generated abelian group $G$. In all that follows, let
$[n+1]=\{1,2,...,n+1\}$ be the first $n+1$ natural numbers, and $p_{i}$ is a
prime for any positive integer $i$. Recall that the pseudoquandle obtained
from a group $G$ is denoted as $P_{G}$.

\begin{proposition}
The pseudoquandle obtained from $%
\mathbb{Z}
/p^{n}%
\mathbb{Z}
$ , with $p$ a prime number, is isomorphic to $[n+1]$ with the operation
$\vartriangle$ defined as $i\vartriangle j=max\{i,j\}$ for all $i,j$ $\in$
$[n+1]$
\end{proposition}

%

\begin{proof}%
As $%
\mathbb{Z}
/p^{n}%
\mathbb{Z}
$ is finite and Abelian, the (normal) subgroups of $%
\mathbb{Z}
/p^{n}%
\mathbb{Z}
$\ are the ones generated by elements whose orders are powers of $p$. We shall
denote these (normal) subgroups by $x_{1},x_{2},...x_{n+1}$\ with
$x_{1}=\{e\}$, $x_{n+1}=%
\mathbb{Z}
/p^{n}%
\mathbb{Z}
$\ so that $x_{i}$ is a subgroup of $x_{j}$\ iff $i<j$\ for all $i,j$\ $\in
$\ $[n+1]$. Also, by definition, $x_{i}\ast x_{j}=\{a+b$ $|$ $a\in x_{i},$
$b\in x_{j}\}$, which is precisely $x_{max\{i,j\}}$. Thus, we need only to
verify that $[n+1]$\ with $\vartriangle$\ is a pseudoquandle. Idempotence is
trivial, while self-distributivity follows from the transitivity of
$max\{-,-\}$.%
\end{proof}%
\medskip

\begin{corollary}
The pseudoquandle obtained from any two finite cyclic groups of the same order
are isomorphic.
\end{corollary}

%

\begin{proof}%
This is immediate from the above proposition.%
\end{proof}%
\medskip

Thus, we see that by creating a pseudoquandle from a finite prime power cyclic
group, we cannot "go back" uniquely to the group since we loose information
about the prime, so that for instance $%
\mathbb{Z}
/p^{n}%
\mathbb{Z}
$ and $%
\mathbb{Z}
/q^{n}%
\mathbb{Z}
$ yield the same pseudoquandle structure for different primes $p$ and $q$.

\textbf{Remark:} We can renumber the elements of the pseudoquandle as $x_{1}=%
\mathbb{Z}
/p^{n}%
\mathbb{Z}
,$ $x_{n+1}=\{e\}$ so that now $x_{i}$ is a subgroup of $x_{j}$ iff $i>j$. By
this labelling, the\ pseudoquandle obtained from $%
\mathbb{Z}
/p^{n}%
\mathbb{Z}
$ is isomorphic to $[n+1]$ with the operation $i\vartriangle j=min\{i,j\}$ for
all $i,j$ $\in$ $[n+1]$. However, $[n+1]$ with $min\{-,-\}$ and with
$max\{-,-\}$ as the product are isomorphic as pseudoquandles so there is no
confusion.\footnote{More generally, if $X$ is an ordered set such that for all
$x,y$ in $X$, $min\{x,y\}$/$max\{x,y\}$ can be defined, one can show that
$(X,\ast)$ with $x\ast y=min/max\{x,y\}$ is a pseudoquandle.}

We can generalize the above to direct sums to obtain:

\begin{proposition}
$P_{G}$ obtained from $G=%
\mathbb{Z}
/p_{1}^{n}%
\mathbb{Z}
\oplus%
\mathbb{Z}
/p_{2}^{m}%
\mathbb{Z}
$ is isomorphic to $[n+1]\oplus\lbrack m+1]$ with binary product given by
$(i_{1},i_{2})\ast(j_{1},j_{2})=(max\{i_{1},i_{2}\},max\{j_{1},j_{2}\})$.
\end{proposition}

%

\begin{proof}%
All the (normal) subgroups of $%
\mathbb{Z}
/p_{1}^{n}%
\mathbb{Z}
\oplus%
\mathbb{Z}
/p_{2}^{m}%
\mathbb{Z}
$ are of the form $(x_{i},y_{j})$ with $x_{i},y_{j}$ (normal) subgroups in $%
\mathbb{Z}
/p_{1}^{n}%
\mathbb{Z}
,%
\mathbb{Z}
/p_{2}^{n}%
\mathbb{Z}
$ respectively with $i,j$ $\in$ $[n+1],[m+1]$ respectively. Thus,
\[
(x_{i_{1}},y_{j_{1}})\ast(x_{i_{2}},y_{j_{2}})=(x_{i_{1}}\ast x_{i_{2}%
},y_{j_{1}}\ast y_{j_{2}})=(x_{max\{i_{1},i_{2}\}},y_{max\{j_{1},j_{2}\}})
\]
with $i_{1},i_{2}$ $\in$ $[n+1],$ $j_{1},j_{2}$ $\in$ $[m+1]$ By the
proposition above, the two factors can each be identified with $(max\{i_{1}%
,i_{2}\},max\{j_{1},j_{2}\})$ bijectively. Hence, we can conclude that the
pseudoquandle obtained from $%
\mathbb{Z}
/p_{1}^{n}%
\mathbb{Z}
\oplus%
\mathbb{Z}
/p_{2}^{m}%
\mathbb{Z}
$ is isomorphic to $[n+1]\oplus\lbrack m+1]$.%
\end{proof}%
\medskip

Next, we consider the pseudoquandle obtained from $G\ =%
\mathbb{Z}
$.

\begin{proposition}
$P_{%
\mathbb{Z}
}$ is isomorphic to $%
\mathbb{Z}
_{+}$ $($as a pseudoquandle$)$ with product $n\ast m=gcd\{n,m\}$.
\end{proposition}

%

\begin{proof}%
This is trivial since for any two subgroups $n%
\mathbb{Z}
,m%
\mathbb{Z}
$ of $%
\mathbb{Z}
,$ we have that
\[
n%
\mathbb{Z}
\ast m%
\mathbb{Z}
=n%
\mathbb{Z}
+m%
\mathbb{Z}
=gcd\{n,m\}%
\mathbb{Z}%
\]
which can be identified with $gcd\{n,m\}$%
\end{proof}%
\medskip

\begin{corollary}
As pseudoquandles, $P_{%
\mathbb{Z}
}\oplus P_{%
\mathbb{Z}
/p^{n}%
\mathbb{Z}
}$ is isomorphic to $%
\mathbb{Z}
_{+}\oplus\lbrack n+1]$
\end{corollary}

%

\begin{proof}%
This follows readily from the above two propositions.%
\end{proof}%
\medskip

Note that the above corollary holds for a finite direct sum of $P_{G}$
terms.\bigskip

We are now ready to state the main result of this section which classifies
$P_{G}$ for any finitely generated abelian group $G$.

\begin{theorem}
\label{1}The pseudoquandle obtained from any finitely generated abelian group
$G$ is isomorphic to $L_{n,r}$\ = $%
\mathbb{Z}
_{+}^{n}\oplus\lbrack m_{1}+1]\oplus\lbrack m_{2}+1]...\oplus\lbrack m_{r}%
+1]$\ with binary product given by$:$%
\[
(x_{1},...x_{n},i_{1},...i_{r})\ast(y_{1},...y_{n},j_{1},...j_{r}%
)=(gcd\{x_{1},y_{1}\},...,gcd\{x_{n},y_{n}\},max\{i_{1},j_{1}\},...,max\{i_{r}%
,j_{r}\})
\]
For all $x_{s},y_{s}$\ $\in$\ $%
\mathbb{Z}
$\ and $i_{t},j_{t}$\ $\epsilon$\ $[m_{t}+1]$ and positive integers
$n,m_{1},m_{2},...,m_{r}$
\end{theorem}

%

\begin{proof}%
We shall prove the theorem by using the fundamental theorem of finitely
generated abelian groups (primary decomposition form). Any finitely generated
abelian group $G$\ is isomorphic to
\[%
\mathbb{Z}
^{n}\oplus%
\mathbb{Z}
/p_{1}^{m_{1}}%
\mathbb{Z}
\oplus%
\mathbb{Z}
/p_{2}^{m_{2}}%
\mathbb{Z}
\oplus...\oplus%
\mathbb{Z}
/p_{r}^{m_{r}}%
\mathbb{Z}%
\]
From the above two propositions and corollary, it is clear that the
pseudoquandle $P_{G}$\ is isomorphic to $%
\mathbb{Z}
_{+}^{n}\oplus\lbrack m_{1}+1]\oplus\lbrack m_{2}+1]...\oplus\lbrack m_{r}%
+1]$\ with the product as given in the theorem.%
\end{proof}%
\medskip

\subsection{The pseudoquandle matrix}

In this section, we will give some characterization of the matrix of $P_{G}$
for arbitrary (finite) groups $G$. The main definitions are similar to those
in the quandle case (see \cite{NH}, \cite{NV}), but we are able to prove more
general results in the case of pseduquandles of the form $P_{G}$.

We recall some basic defintions pertaining to the quandle matrix before
continuing with our discussion. The quandle/rack/pseudoquandle matrix is
simply the multiplication table of the corresponding object written in matrix form:

\begin{definition}
Let $X$ be a rack $($resp. quandle, pseudoquandle$)$ so that $X=\{x_{1}%
,x_{2},...,x_{n}\}$ for some $n$. Define the $X$-matrix $M_{X}$ as the
following$:$%
\[
M_{X}=%
\begin{pmatrix}
x_{1}\ast x_{1} & x_{1}\ast x_{2} & ... & ... & x_{1}\ast x_{n}\\
x_{2}\ast x_{1} & x_{2}\ast x_{2} &  &  & x_{2}\ast x_{n}\\
... &  & ... &  & ...\\
... &  &  & ... & ...\\
x_{n}\ast x_{1} & x_{n}\ast x_{1} & ... & ... & x_{n}\ast x_{n}%
\end{pmatrix}
\]

\end{definition}

As in the case of the quandle matrix, we shall use only the subscript indices
to denote the elements of the pseudoquandle matrix to stay consistent with the
notation used in literature. Thus $M_{P_{G}}$ will be an integral matrix, so
that the $(i,j)$ element is the subscript of the element $x_{i}\ast x_{j}$. We
will always consider the $(1,1)$ element of the matrix $M_{P_{G}}$ to be $1$,
corresponding to the trivial subgroup.

\begin{proposition}
Let $G$\ be any $($finite$)$ group. Let $P_{G}$ be pseudoquandle obtained from
$G$ and $M_{P_{G}}$ be the corresponding integral pseudoquandle matrix. Then,
$M_{P_{G}}$ can be diagonalized.
\end{proposition}

%

\begin{proof}%
Since $P_{G}$\ is commutative for every group $G$, we can see that $M_{P_{G}}%
$\ is a symmetric matrix and hence can be diagonlized by an orthogonal matrix
by the spectral theorem.%
\end{proof}%
\medskip

\begin{corollary}
$M_{P_{G}}$ is of the form
\[%
\begin{pmatrix}
1 & 2\\
2 & 2
\end{pmatrix}
\]
iff $G$ is simple.
\end{corollary}

%

\begin{proof}%
If $G$ is simple, $M_{P_{G}}$ is a $2\times2$ matrix three of whose elements
are $G$ and the other the trivial group in the $(1,1)$ position proving one
direction. If $M_{P_{G}}$ is of the form given in the proposition, then for
each pair of normal subgroup in $P_{G}$, $x_{i},x_{j}$ the product of
$x_{i},x_{j}$ must some fixed normal subgroup $x_{k}$, taking $x_{i}=\{e\},$
$x_{j}=G$ we see $x_{k}$ must be $G$ for each pair $(i,j)$, which is possible
only if $G$\ is simple.%
\[
M_{P_{G}}=%
\begin{pmatrix}
\{e\} & G\\
G & G
\end{pmatrix}
\leftrightarrows%
\begin{pmatrix}
1 & 2\\
2 & 2
\end{pmatrix}
\]%
\end{proof}%
\medskip

As in the case of quandle matrices (see \cite{NH}), the trace of $M_{P_{G}}$is
$n(n+1)/2$ with $n=|P_{G}|$.

In \cite{LR} it was noted that analogous ideas may be developed for
quandles/racks via permutations of $[n]$ (so that one studies the underlying
operation as an element of $S_{n}$). This does not yield much fruit in the
pseudoquandle case since these do not satisfy the bijection property that
quandles/racks do and thus, the map $\ast_{i}:[n]\rightarrow\lbrack n]$
defined, as usual, by $\ast_{i}(j)=j\ast i$ may not correspond to a
permutation of $[n]$ in the case of pseudoquandles.

\section{Commutative pseudoquandles and kernels}

In this section, we define the main algebraic structure of our study, the
kernel $ker(p)$ of an element $p$ in a commutative pseudoquandle $P$. We will
prove some properties of $ker(p)$ which motivate our two main results, namely,
establishing a bound on the cardinality of certain psedoquandles via kernels
and obtaining a "class equation" for pseudoquandles satisfying an acending
chain criterion (defined later in this section).

\begin{definition}
Given a pseudoquandle $P$, the \textbf{kernel} of an element $p\in P$ is
defined as follows$:$%
\[
ker(p)=\{q\in P|\text{ }p\ast q=q\ast p=p\}
\]
The \textbf{cokernel} $coker(p)$ of $p\in P,$ is defined as$:$%
\[
coker(p)=P-ker(p)=\{q\in P|\text{ }p\ast q=q\ast p\neq p\}
\]

\end{definition}

In the context of quandle polynomials (see \cite{N-2}), the cardinality of
kernels were used to define the polynomial invariants of quandles but the
underlying structure was not analyzed. The results presented there can be
generalized via the notion of kernels. We will have more to say in the case of
(commutative) pseudoquandles.

In all that follows, let $P$ be a pseudoquandle. A subset $R$ of $P$ is a
sub-pseudoquandle if $R$ is a pseudoquandle in its own right. Clearly, any
subset closed with respect to the operation in $P$ is a sub-pseudoquandle
since self-distributivity and idempotence is obtained from the structure of
$P$

\begin{proposition}
$ker(p)$ is a sub-pseudoquandle for each $p\in P$.
\end{proposition}

%

\begin{proof}%
If $x,y\in ker(p)$, then $p\ast(x\ast y)=(p\ast x)\ast(p\ast y)=p\ast p=p$. So
that $x\ast y\in ker(p)$. Thus, $ker(p)$ is closed and a sub-pseudoquandle for
each $p\in ker(p)$.%
\end{proof}%
\medskip

Note that $p\in ker(p)$ for all $p\in P$ and also that $P=ker(p)\sqcup
coker(p)$ (disjoint union). Also, $coker(p)$ need not be a sub-pseudoquandle.
For example, consider the commutative pseudoquandle $P=\{x_{1},x_{2},x_{3}\}$
with the relation $x_{i}\ast x_{j}=x_{k}$, $i\neq j\neq k$ and $x_{i}\ast
x_{i}=x_{i}$ (this is actually the dihedral quandle on $%
\mathbb{Z}
/3%
\mathbb{Z}
$). Here, $ker(x_{1})=\{x_{1}\}$. But $coker(x_{1})=\{x_{2},x_{3}\}$ is not a
sub-pseudoquandle since $x_{2}\ast x_{3}=x_{1}\notin coker(x_{1})$. However,
there is a sufficient condition to ensure that $coker(p)$ is a sub-pseudoquandle:

\begin{proposition}
Let $P=\{p_{1},p_{2},...,p_{n}\}$. If $ker(p_{1})\subseteq ker(p_{2}%
)\subseteq...\subseteq ker(p_{n})$ then, $coker(p_{i})$ is a sub-pseudoquandle
for $i=1,2,...n$.
\end{proposition}

%

\begin{proof}%
Since $p_{1}\in ker(p_{1})\subseteq ker(p_{2})$, $p_{1}\in ker(p_{2})$.
Likewise, $p_{1}$ and $p_{2}$ $\in$ $ker(p_{3})$ etc. so that
\[
p_{1},p_{2},...,p_{n}\in ker(p_{n})
\]
But also, $ker(p_{n})\subseteq P$, which shows that $P=ker(p_{n})$. Let us
assume that $x,y\in coker(p_{i})$. This implies that $p_{i}\in ker(x)$ and
$p_{i}\in ker(y)$ by the ascending kernel assumption. Thus,
\[
p_{i}\ast(x\ast y)=(p_{i}\ast x)\ast(p_{i}\ast y)=x\ast y\neq p_{i}%
\]
so that $x\ast y\in coker(p_{i})$, proving the proposition.%
\end{proof}%
\medskip

\textbf{Remark:} The condition $ker(p_{1})\subseteq ker(p_{2})\subseteq
...\subseteq ker(p_{n})$ is satisfied, for example, by $P_{G}$ with $G=%
\mathbb{Z}
/p^{n-1}%
\mathbb{Z}
$ for positive integers $n$. We will refer to this as the \textbf{ascending
chain criterion}. Note that this implies that $P=ker(p_{n})$.

Next, we show what happens in the intersection of two kernels. It is
interesting to note that the intersection of kernels may not always be a
kernel of another element of the pseudoquandle.

\begin{proposition}
For $p,q\in P$, $ker(p)\cap ker(q)\subseteq ker(p\ast q)$
\end{proposition}

%

\begin{proof}%
If $x\in ker(p)\cap ker(q)$, $x\ast p=p$ and $x\ast q=q$. Thus, $x\ast(p\ast
q)=(x\ast p)\ast(x\ast q)=p\ast q$ so that $x\in ker(p\ast q)$.%
\end{proof}%
\medskip

We will now prove a few results analyzing the structure of finite, commutative
quandles via kernels. The discussion will motivate the class equation for
pseudoquandles satisfying the ascending chain criterion.

\begin{proposition}
For $p,q\in P$ if $p\in ker(q)$ then $ker(p)\subseteq ker(q)$
\end{proposition}

%

\begin{proof}%
Let $x\in ker(p)$. We need to show that $x\ast q=q$. Now, $p\ast q=q$, and
$x\ast p=p$. We also have the following:%
\[
x\ast q=x\ast(p\ast q)=(x\ast p)\ast(x\ast q)=p\ast(x\ast q)
\]
But,
\[
p\ast(x\ast q)=(p\ast x)\ast(p\ast q)=p\ast q=q
\]
so that $x\ast q=q$.%
\end{proof}%
\medskip

Analogous to any other algebraic structure, we can defined the product of two
subsets of a quandle/rack/pseudoquandle etc. with respect to the operation
binary operation inherited from the parent structure.

\begin{definition}
For any subsets $A,B$ of $P$, define
\[
A\ast B=\{a\ast b|\text{ }a\in A,\text{ }b\in B\}
\]
and
\[
A^{2}=A\ast A
\]

\end{definition}

We have that the $ker(p)$ is idempotent as a set with the above definitions:

\begin{proposition}
For all $p\in P$, $[ker(p)]^{2}=ker(p)$
\end{proposition}

%

\begin{proof}%
Let $ker(p)=\{p_{1},p_{2},...,p_{k}\}$. Thus, $ker(p)\ast ker(p)=\{p_{i}\ast
p_{j}|$ $i,j=1,2,...,k\}$. Since $p_{i}\ast p_{i}=p_{i}$, $ker(p)\subseteq
\lbrack ker(p)]^{2}$. Also $ker(p)$ is a sub-pseudoquandle, so that $p_{i}\ast
p_{j}\in ker(p)$ for all $p_{i},p_{j}$ $\in ker(p)$, showing that
$[ker(p)]^{2}\subseteq ker(p)$. Thus, $[ker(p)]^{2}=ker(p)$%
\end{proof}%
\medskip

If $\{p\}$ denotes the single element sub-pseudoquandle, it is easy to see
that $\{p\}\ast ker(p)=\{p\}$. Infact, more is true:

\begin{proposition}
For any $p,q\in P$, $\{q\}\ast ker(p)$ is a sub-pseudoquandle
\end{proposition}

%

\begin{proof}%
Let $x,y\in\{q\}\ast ker(p)$. Thus, $x=q\ast p_{1}$ and $y=q\ast p_{2}$ where
$p_{1}$ $and$ $p_{2}\in ker(p)$. Clearly, $x\ast y=q\ast(p_{1}\ast p_{2})$
which is an element of $\{q\}\ast ker(p)$. Like before, idempotence and
self-distributivity follow from multiplication in $P$.%
\end{proof}%
\medskip

The next lemma will be required for establishing a lower bound on the
cardinality of a pseudoquandle whose kernels and cokernels satisfy certain
intersection criteria.

\begin{lemma}
If $p,q\in P$, and $ker(p)\cap ker(q)=\phi$, then $ker(q)\subseteq coker(p)$
and $ker(p)\subseteq coker(q)$\newline
\end{lemma}

%

\begin{proof}%
Since $P=ker(p)\sqcup coker(p),$%
\[
P\cap ker(q)=(ker(p)\sqcup coker(p))\cap ker(q)=(ker(p)\cap ker(q))\sqcup
(coker(p)\cap ker(q))
\]
but $ker(p)\cap ker(q)=\phi$ by our assumption. Thus $ker(q)=coker(p)\cap
ker(q)$ showing that
\[
ker(q)\subseteq coker(p)
\]
The other inclusion is entirely similar.%
\end{proof}%
\medskip

A simple application of the above proposition gives us the lower bound on the
cardinality of $P$, when $|ker(p)|=|ker(q)|=k$.

\begin{corollary}
If $p,q\in P$, and $ker(p)\cap ker(q)=\phi$ with $|ker(p)|=|ker(q)|=k$, then
$|P|\geq2k$
\end{corollary}

%

\begin{proof}%
$P=ker(p)\sqcup coker(p)$. So $|P|=|ker(p)|+|coker(p)|$. But by the above
corollary, $coker(p)$ contains $ker(q),$ so that $|coker(p)|\geq|ker(q)|$.
Hence, $|P|\geq|ker(p)|+|ker(q)|=2k$.%
\end{proof}%
\medskip

Let us define
\[
P^{ker}=\{ker(p)|\text{ }p\in P\}
\]
We will now show that the map $\varphi:P\rightarrow P^{ker}$ defined by
$\varphi(p)=ker(p)$ is bijective.

\begin{proposition}
The map $\varphi$ defined above is bijective
\end{proposition}

%

\begin{proof}%
We will first show that $p=q\iff ker(p)=ker(q)$. If $p=q$, it is clear that
$ker(p)=ker(q)$. For the other direction, $ker(p)=ker(q)$ implies that $q\in
ker(p)$ and $p\in ker(q)$, so $p=q\ast p=p\ast q=q$. Hence, $\varphi$ is
injective. Also, $\varphi$ is clearly surjective hence is a bijection.%
\end{proof}%
\medskip

We now come to the other main result of this section. Before this, we will
define a few terms used in the result.

\begin{definition}
Let $p\in P$ and $q\in ker(p)$. The \textbf{relative cokernel} of $q$ in $p$
is $ker(p)-ker(q)$ written as $coker(q:p)$. If $q\notin ker(p)$, we define
$coker(q:p)=\phi$. Thus
\[
coker(q:p)=\{s\in ker(p)|s\ast q\neq q\}=ker(p)-ker(q)
\]
The cardinality of $coker(q:p)$ is defined to be the \textbf{index} of $q$ in
$p$.
\end{definition}

The following theorem is a pseudoquandle version of Lagrange's theorem in
group theory.

\begin{theorem}
Let $p\in P$ and $q\in ker(p)$. Then, $|P|=|ker(q)|+|coker(p)|+|coker(q:p)|$.
\end{theorem}

%

\begin{proof}%
$P=ker(p)\sqcup coker(p)$ so that%
\[
P-ker(q)=(coker(p)\sqcup ker(p))-ker(q)
\]%
\[
\implies coker(q)=coker(p)\sqcup(ker(p)-ker(q))
\]%
\[
\implies coker(q)=coker(p)\sqcup coker(q:p)
\]
But also, $P=ker(q)\sqcup coker(q)=ker(q)\sqcup(coker(p)\sqcup coker(q:p))$.
Taking cardinalities,
\[
|P|=|ker(q)|+|coker(p)|+|coker(q:p)|
\]%
\end{proof}%
\medskip

The following is an application of the above ideas and is a "class equation"
of pseudoquandles satisfying the acending chain criterion. It is similar to
the class equation for finite groups in terms of (the cardinalities of)
conjugacy classes.

\begin{theorem}
Let $P=\{p_{1},p_{2},...,p_{n}\}$ with $ker(p_{1})\subseteq ker(p_{2}%
)\subseteq...\subseteq ker(p_{n})$. Then
\[
|P|=|ker(p_{1})|+%
{\displaystyle\sum\limits_{k=1}^{n-1}}
|coker(p_{k}:p_{k+1})|
\]

\end{theorem}

%

\begin{proof}%
It was remarked above that that $ker(p_{1})\subseteq ker(p_{2})\subseteq
...\subseteq ker(p_{n})$ implies $P=ker(p_{n})$. Now, we have the following
manipulation:
\[
ker(p_{n})=(ker(p_{n})-ker(p_{n-1}))\sqcup...\sqcup(ker(p_{2})-ker(p_{1}%
))\sqcup(ker(p_{1}))
\]
Taking cardinalities on both sides, we see that $|P|=|ker(p_{n})|=|ker(p_{1}%
)|+%
{\displaystyle\sum\limits_{k=1}^{n-1}}
|coker(p_{k}:p_{k+1})|$.%
\end{proof}%
\medskip

We close the section with the following observation regarding the behaviour of
kernels of a pseudoquandle under a homomorphism:

\begin{proposition}
Let $\theta:P\rightarrow Q$ be a homomorphism of pseudoquandles. Then,
$\theta(ker(p))\subseteq ker(\theta(p))$ for each $p\in P$. There is equality
of sets if $\theta$ is an isomorphism.
\end{proposition}

%

\begin{proof}%
Let $x\in ker(p)$. Then $\theta(x)\ast\theta(p)=\theta(x\ast p)=\theta(p)$.
Hence, $\theta(x)\in ker(\theta(p))$. If $\theta$ is an isomorphism, and
$\theta(x)\in ker(\theta(p))$ then clearly $\theta(x)\ast\theta(p)=\theta
(p)=\theta(x\ast p)$ which implies $x\ast p=p$, thus they two sets are
equivalent.%
\end{proof}%
\medskip

The restriction of an isomorphism between two pseudoquandles to kernels
results in an induced isomorphism. If $P$ and $Q$ are isomorphic, say via
$\Phi$, the restriction of $\Phi$ to each $ker(p_{i})$ induces an isomorphism
$\Phi|_{ker(p_{i})}$ between the kernels $ker(p_{i})$ and $ker(q_{i})$ for
each $i=1,2,...n$. Clearly, for $p\in ker(p_{i})\cap ker(p_{j})$, $\varphi
_{i}(p)=\varphi_{j}(p)$.

The converse to the above may not always true. However, it is easily seen to
be true in the case of pseudoquandles satisfying the ascending chain
criterion. In that case, we can characterize (upto isomorphism) a
pseudoquandle from its constituent kernels.

\end{document}